\documentstyle[11pt]{article}
\newcommand{\ds}{\displaystyle}
\newtheorem{theorem}{Theorem} [section]
\newtheorem{lemma}[theorem]{Lemma}
\newtheorem{corollary}[theorem]{Corollary}

\newtheorem{proposition}[theorem]{Proposition}
\newtheorem{example}[theorem]{Example}

\newcommand{\norm}[1]{||} 

\newcommand{\dd}{{\mathrm d}}  

\newcommand{\RR}{{\bf R}}  

\newcommand{\Ll}{{\mathcal L}}  
\newcommand{\Hh}{{\mathcal H}}  
\newcommand{\Vv}{{\mathcal V}}  

\newcommand{\ii}{{\rm i}}  

\newcommand{\ra}{\rightarrow}

\newcommand{\pa}{\partial}

\renewcommand{\phi}{\varphi}
\newcommand{\la}{\lambda}

\newcommand{\na}{\nabla}
\newcommand{\al}{\alpha}
\newcommand{\be}{\beta}
\newcommand{\ga}{\gamma}

\newcommand{\si}{\sigma}

\newcommand{\om}{\omega}
\newcommand{\Om}{\Omega}

\newcommand{\wt}{\widetilde} 
\newcommand{\ov}{\overline} 

\renewcommand{\Re}{{\rm Re}\,} 
\renewcommand{\div}{\mbox{\rm div}\,} 
\newcommand{\grad}{\mbox{\rm grad}\,}
\newcommand{\Tr}{\mbox{\rm Tr}\,}
\newcommand{\Ric}{\mbox{\rm Ricci}\,}

\begin{document}
\title{Three-dimensional Ricci solitons which project to surfaces}
\author{Paul Baird and Laurent Danielo}
\date{}
\maketitle{\small \begin{center}D\'epartement de Math\'ematiques, Universit\'e de Bretagne Occidentale \\6 Avenue Le Gorgeu, B.P. 452, 29285 Brest Cedex, France  \\E-mail:Paul.Baird@univ-brest.fr \ \  Laurent.Danielo@univ-brest.fr \\{\it 2000 Mathematical Subject Classification:}  57M50, 35Q51 \\{\it Key words:}  Ricci soliton, semi-conformal map\end{center}  }

\begin{abstract}We study $3$-dimensional Ricci solitons which project via a semi-conformal mapping to a surface.  We reformulate the equations in terms of parameters of the map; this enables us to give an ansatz for constructing solitons in terms of data on the surface.  A complete description of the soliton structures on all the $3$-dimensional geometries is given, in particular, non-gradient solitons are found on Nil and Sol. 
\end{abstract}

\section{Introduction}Fixed points of the Ricci flow:$$\frac{\pa g}{\pa t} = - 2\Ric (g)$$on a manifold $(M, g)$ are called \emph{Ricci solitons}.  It is usual to look for fixed points up to \emph{diffeomorphism} of $M$ and \emph{scaling} of $g$, whence the equations for a Ricci soliton become:
\begin{equation} \label{eq:soliton}
- 2 \Ric (g) = \Ll_Eg + 2Ag
\end{equation}
where $E$ is a vector field on $M$, $\Ll_Eg$ denotes Lie derivation of $g$ with respect to $E$ and $A$ is a constant.  The soliton is called \emph{shrinking, stationary} or \emph{expanding} according as $A<0, =0, >0$, respectively.  It is of \emph{gradient type} if $E = \grad F$ for some function $F$, in which case $\Ll_Eg = 2 \na \dd F$.  

Apart from their interest as fixed points of the Ricci flow, see \cite{Ha}, and \cite{Kn} for an overview, they are interesting geometric objects in their own right.   Warped product solutions have been constructed by Bryant and Ivey \cite{Iv-2}. K\"ahler solitons in even dimension are studied by Koiso \cite{Ko}, Cao \cite{Ca-1, Ca-2}, Feldman, Ilmanen and Knopf \cite{Il-Kn} and Bryant \cite{Br-2}.  A construction using a doubly warped product metric in higher dimensions is given by Gastel and Krong \cite{Ga-Kr}.  We note that all of the above constructions are of gradient type.  By a result of Ivey \cite{Iv-1}, in dimension $3$ the only compact examples are of constant curvature.    Our aim in this article is to construct $3$-dimensional solitons which admit a semi-conformal projection onto a surface.  Some of our examples are not of gradient type.  The existence of such solitons has important consequences for the  stability of the Ricci flow about non-Ricci flat metrics, which is studied in the Ricci-flat case in \cite{Gu-Is-Kn, Se}.

A Lipschitz map $\phi : (M^m, g) \ra (N^n, h)$ between Riemannian manifolds is said to be \emph{semi-conformal} if, at each point $x\in M$ where $\phi$ is differentiable (dense by Radmacher's Theorem), the derivative $\dd \phi_x : T_xM\ra T_{\phi (x)}N$ is either the zero map or is conformal  and surjective on the compliment of $\ker \dd\phi_x$ (called the \emph{horizontal distribution}).  Thus, there exists a number $\la (x)$ (defined almost everywhere), called the \emph{dilation}, such that $\la (x)^2 g(X,Y)$ $ = $ $\phi^*h(X,Y)$, for all $X,Y\in (\ker \dd\phi_x)^{\perp}$.  If $\phi$ is of class $C^1$, then we have a useful characterisation in local coordinates, given by
$$
g^{ij}\phi_i^{\al}\phi_j^{\be} = \la^2 h^{\al\be}\,,
$$ 
 where $(x^i), (y^{\al})$ are coordinates on $M, N$, respectively and $\phi^{\al}_i = \pa (y^{\al}\circ \phi ) / \pa x^i$ (here and throughout we sum over repeated indices).  The  fibres of a smooth submersive semi-conformal map determine a conformal foliation, see \cite{Va} and conversely, with respect to a local foliated chart, we may put a conformal structure on the leaf space with respect to which the projection is a semi-conformal map.  We then have the identity:
$$
\left(\Ll_Ug\right)(X,Y) = - 2 U(\ln \la )\,g(X,Y),
$$
for $U$ tangent and $X,Y$ orthogonal to the foliation.

The question of which $3$-manifolds admit a conformal foliation (equivalently, local semi-conformal maps to surfaces) remains a delicate one.  It is equivalent, locally, to finding a function which admits a `conjugate'; in \cite{Ba-Ea} it is shown that such functions are characterized as satisfying a 2nd order differential inequality and a 3rd order differential equation.  We shall see that we can calculate the Ricci curvature in terms of geometric quantities associated to the projection, for example the mean curvature of the fibres, the integrability tensor of the orthogonal distribution and the dilation.  

In general, the equations are still difficult to study and in this article we concentrate on the case when objects are \emph{basic}, that is they are defined in terms of data on the surface.  This leads to the following ansatz for the construction of solitons.

\begin{theorem} \label{th:main}  Let $(N,h)$ be a Riemannian surface. 
Let $\ov{\psi}, \ov{\la}, \ov{\rho}, \ov{\nu} :N\ra \RR$ be functions with $\ov{\psi}\geq 0$, satisfying the system of equations
\begin{equation}\label{eq:main}
\left\{ \begin{array}{rrcl}
{\rm (i)} & K^N +  \frac{1}{2}\left( \Delta^N\ln (\ov{\la}^2\ov{\nu}) - |\grad \ln\ov{\rho}|^2\right) +  \frac{A-\ov{\psi}}{\ov{\la}^2} & = & 0 \\
{\rm (ii)(a)} &  \Delta^N\ln (\ov{\rho}^2\ov{\nu}\ov{\psi}^{-1/2}) +   |\grad \ln\ov{\rho}|^2 
  - \frac{1}{4}|\grad \ln\ov{\psi}|^2  \qquad \qquad \\ 
  &  + \frac{1}{2}h(\grad\ln\ov{\psi} ,\grad\ln \ov{\nu})+ \frac{\ov{\psi} + A}{\ov{\la}^2} & = & 0 \\
{\rm (ii)(b)}  & \ov{\la}^2 \left\{ \Delta^N\ln\ov{\rho} - h(\grad \ln\ov{\rho}, \grad\ln \ov{\nu}) + \frac{\ov{\psi} + A}{\ov{\la}^2}\right\} &  = & {\rm const.} \\
{\rm (u)} &  \quad \nabla \dd\ln\ov{\nu} + 2 \dd \ln\ov{\la}\odot \dd\ln\ov{\nu} - (\dd\ln\ov{\rho} )^2  =  \al h \ \ ({\rm some} \ \al : N\ra \RR )\,, & &   
\end{array} \right.
\end{equation}
where we take {\rm (ii)(a)} to be vacuous whenever $\ov{\psi} \equiv 0$ and
where $K^N$ is the Gaussian curvature of $N$.  Now set $M = N\times (-\delta , \delta )$ for some $\delta >0$ and let $\phi : M\ra N$ be the canonical projection.  Let $\ov{\si} = \sqrt{2\ov{\psi}}/\ov{\la}^2$ and set $\ov{\Omega} = \ov{\si} \mu^N$, where $\mu^N$ denotes the volume form on $N$.  Write $\wt{\Omega} = \phi^*\ov{\Omega}$, $\rho = \ov{\rho}\circ \phi$ and $\la = \ov{\la}\circ \phi$.  Let the $1$-form $\theta$ be a solution to the exterior differential equation 
\begin{equation} \label{eq:ext}
\dd\theta + \dd\ln\rho \wedge \theta = \wt{\Omega}
\end{equation}
which is everywhere non-vanishing on $\ker \dd\phi$.
 Write $g = \frac{\phi^*h}{\la^2} + \theta^2$.  Then $(M^3, g)$ is a Ricci soliton.
\end{theorem}

The function $\ov{\psi}$ measures the non-integrablility of the orthogonal compliment of $\ker \dd\phi$, in particular if $\ov{\psi}\equiv 0$, then this distribution is integrable.  The function $\ln\rho$ is the \emph{potential} for the mean curvature $\mu$ of the fibres of $\phi$; that is, $\mu = \grad \ln\rho$. 
A special case of the theorem is when $\ov{\rho}$ is constant, so the fibres are minimal.   In this case we can give a complete description of the solutions of (\ref{eq:main}) in terms of holomorphic data on the surface $N$, which leads to the following conclusion.

\begin{corollary} \label{cor:main}
Let $(M^3, g)$ be a soliton derived from the ansatz given by {\rm Theorem \ref{th:main}} with $\rho$ constant, then $\ov{\psi}=C$ is constant and either

{\rm (i)} $C= 0$, in which case $M$ is locally isometric to a Riemannian product $N^2 \times \RR$, where $N^2$ is $2$-dimensional gradient Ricci soliton, or 

{\rm (ii)}  $C\neq 0$, in which case $M$ either has constant curvature or is locally isometric to the geometry {\rm Nil}.
\end{corollary}

We will give an explicit description of the soliton structure on Nil, which exists and is unique by \cite{La}.
 On the other hand, to identify the geometry Sol as a soliton (Example \ref{ex:sol}) we will require the more general set-up of Theorem \ref{th:main}.   This is because Sol, even locally, does not admit a semi-conformal map to a surface with geodesic fibres \cite{Ba-Wo}.   For both of these geometries the soliton structure is not of gradient type.  Finally, we show that $\wt{{\rm SL}_2(\RR )}$ does not admit any soliton structure whatsoever.  To complete the picture, we show uniqueness of the soliton structures on the other $3$-dimensional geometries, except for $\RR^3$ where the Gaussian soliton appears.

\medskip

\noindent {\bf Remarks}:  1.  The constant $A$ is the same constant that occurs in (\ref{eq:soliton}) and so its sign determines whether the soliton is shrinking, stationary or expanding.  To be consistent with notation, functions, forms and vector fields on the codomain $N$ will have a `bar' over them. 

2.  The system (\ref{eq:main}) depends only on $\grad \ln \ov{\nu}$ and $\grad \ln \ov{\rho}$.  In its most general form we may replace $\grad \ln \ov{\nu}$ by $\ov{Y}$, where $\ov{Y}$ is a vector field on $N$.  For example, $\Delta^N\ln\ov{\nu}$ would then become $\div \ov{Y}^{\flat}$, where $\ov{Y}^{\flat} = h(\ov{Y}, \ \cdot \ )$ is the dual of $\ov{Y}$.  Our hypothesis is  not the same as supposing the soliton is of gradient type (see also Section \ref{sec:killing}).

3.  We will refer to the condition (u) in the above theorem and corollary as the `umbilicity condition'; the reason for this will become apparent in Section \ref{sec:killing}.

4.  The system (\ref{eq:main}) is invariant under conformal changes in the metric $h$.  This is a natural consequence of our construction.  Indeed, if we replace $h$ by $e^{2u}h$, then $\ov{\la}$ is replaced by $e^u\ov{\la}$, $K^N$ by $e^{-2u}(-\Delta^Nu+K^N)$ and the Laplacian $\Delta^N$ by $e^{-2u}\Delta^N$.

5.  Locally, on a topologically trivial domain, it is always possible to solve (\ref{eq:ext}).  In fact we can write it in the equivalent form: $\dd (\rho \theta ) = \rho \wt{\Om}$.  In order to solve this for $\theta$ we require that $\dd (\rho \wt{\Om})=0$.  But $\rho \wt{\Om} = \phi^*(\ov{\rho}\ov{\Om})$, which is clearly closed.

6.  The condition that a soliton given by Theorem \ref{th:main} has constant curvature is given in Proposition \ref{prop:cc}.
  
7.  A semi-conformal map which is \emph{harmonic} is called a harmonic morphism (see \cite{Fu, Ba-Wo}).  If the mapping takes values in a surface then this is equivalent to the fibres being minimal.  The second named author has constructed Einstein metrics in dimension $4$ from harmonic morphisms with $1$-dimensional fibres \cite{Da}.  On the other hand, any harmonic morphism with $1$-dimensional fibres from an Einstein manifold of dimension $\geq 5$ is of warped product type or of Killing type (the fibres are the integral curves of a Killing vector field).  This was shown for manifolds of constant curvature by Bryant \cite{Br-1} and for more general Einstein manifolds by Pantilie and Wood \cite{Pa-Wo}.  In dimension $4$ one other type can occur \cite{Pa}.

\medskip

The first named author thanks the Australian Research Council for funding a visit to the School of Pure Mathematics at the University of Adelaide in 2003, during which time much of this work was done; he is also grateful to the school for their hospitality.  We also thank Dan Knopf for many helpful comments.

\section{The fundamental equations} \label{sec:warped}

Let $\phi : (M^3, g) \ra (N^2, h)$ be a ($C^{\infty}$) semi-conformal submersion with dilation $\la$.  We will use the following notation:  Let $\Vv$ and $\Hh$ denote the \emph{vertical} and \emph{horizontal} distributions, respectively.  We will use the same letters to denote orthogonal projection onto the respective distributions.  Let $U$ denote the \emph{unit vertical vector field}, i.e. $\dd\phi (U) = 0$ and $g(U,U) = 1$.  If $M$ is oriented, there are two choices for $U$;  we will select one of these - the equations are invariant of this choice.  If $M$ is not oriented, then we work locally with a choice of $U$.  Let $\theta = U^{\flat} = g(U, \ \cdot \ ) = g\rfloor U$ denote its dual and write $\Omega = \dd\theta$; then the integrability tensor $I(X,Y) = \Vv [X,Y],\, X,Y \in \Hh$ is related by the formula $\Om (X,Y) = - g(I(X,Y), U)$.  We will use the norm: $||\Om ||^2 = \sum_{a,b}|\Om (e_a, e_b)|^2$, where $\{ e_a\}$ is an orthonormal frame, similarly for $||I||^2$.  Write $\dd^{\Vv}f$ for the vertical component of the exterior differential of a function $f$, i.e. $\dd^{\Vv}f = U(f)\theta$.  Set $\om \odot \eta = \frac{1}{2}(\om \otimes \eta + \eta \otimes \om )$ for the symmetric product of two $1$-forms.   Finally, our sign convention for the Laplacian is such that $\Delta f = f^{\prime\prime}$ for a function of a single variable. 
\begin{proposition} \label{pro:Ricci}  Let $\phi : (M^3, g) \ra (N^2, h)$ be a ($C^{\infty}$) semi-conformal submersion with dilation $\la$.   Then, in the notations defined above,  the Ricci tensor of $(M, g)$ is given by the formula:
\begin{eqnarray} \label{eq:Ricci}
\Ric (g)  & = & \left\{ \la^2 K^N + \Delta \ln \la + \mu (\ln \la )\right\} (g - \theta^2) - \frac{1}{4} ||I||^2 g   + \frac{1}{2}\Ll_{\mu} g \nonumber \\
  &  &  - (\mu^{\flat} + \dd^{\Vv}\ln \la )^2 - (\dd^{\Vv}\ln\la )^2 + 2 \dd (U(\ln\la ))\odot \theta + \dd^*\Om \odot \theta ,  
\end{eqnarray}   
where $\mu$ denotes the mean-curvature vector field of the fibres of $\phi$ and $K^N$ is the Gauss curvature of $N$.
\end{proposition}

The above formula will be deduced by evaluating the Ricci tensor on different combinations of horizontal and vertical vectors.
In what follows, we let $\{ e_i\} = \{ e_a, U\}$ denote a local orthonormal frame field, where the index $i$ ranges over $1,2,3$ and the index $a$ over $1,2$.  The following lemma is a straightforward calculation.

\begin{lemma} \label{lem:app-1} Let $\wt{\Om}$ be the $2$-form defined by $\wt{\Om}(E,F) = - g(U, [\Hh E, \Hh F])$, for all vectors $E,F$. Then
$$
\Om = \wt{\Om} - \mu^{\flat}\wedge \theta\,.
$$
\end{lemma}

We will call $\wt{\Om}$ the \emph{integrablility $2$-form associated to} $\phi$.
The following lemma is useful and is to be found in \cite{Ba}.

\begin{lemma} \label{lem:Hmc} The mean curvature $\ds \frac{1}{2}g(U, \na_{e_a}e_a)$ of the horizontal distribution $\Hh$ is given by $U(\ln\la )$.
\end{lemma}

\noindent \emph{Proof}:  We note that the quantity $g(U, \na_{e_a}e_a)$ is independent of the horizontal frame $\{ e_a\}$.  Suppose that, for each $a = 1,2$, the vector field $e_a$ is the normalised lift of a unit vector field $\ov{X}_a$ on a domain of $N$:  $\dd\phi (e_a) = \la \ov{X}_a\circ \phi$.  By the symmetry of the second fundamental form $\na \dd \phi$, we have
$$
- \dd \phi (\na_{e_a}U) = - \dd \phi (\na_Ue_a) + \na_U^{\phi^{-1}TN}\dd\phi (e_a)\,,
$$
so that $\Hh [U, e_a] = U(\ln\la )e_a$.  
Thus
$$
g(U, \na_{e_a}e_a)  =  g([U, e_a], e_a) 
  =  g(U(\ln\la )e_a, e_a) = 2U(\ln\la )\,.
$$
 \hfill q.e.d. 

 For horizontal vectors $X,Y$, let $A_XY = \Vv \na_XY = \frac{1}{2}\Vv [X,Y] + g(X,Y) \Vv \grad \ln\la $ be one of the fundamental tensors associated to a submersion.  The second fundamental tensor $B$ is defined by $B_UV = \Hh \na_UV$ for vertical vectors $U,V$  \cite{Ba-Wo, ON}.  We write $A^*, B^*$ for their respective adjoints.

The formula for the Ricci tensor evaluated on vertical vectors is as follows.
            
\begin{lemma} \label{lem:vert} $\Ric (U,U) = 2 U(U(\ln\la )) - 2 U(\ln\la )^2 - \frac{1}{4}||I||^2 + \frac{1}{2}(\Ll_{\mu}g) (U,U) + \dd^*\Om (U)$.
\end{lemma}
            
\noindent \emph{Proof}:  First note that 
$$
\Ric (U,U) = \sum_a g(R(U, e_a)e_a, U)) = \sum_a K (e_a\wedge U)\,,
$$
where $K(e_a\wedge U)$ is the sectional curvature of the plane spanned by $e_a$ and $U$.  By Proposition 11.2.2 of \cite{Ba-Wo}, this latter term equals:
\begin{eqnarray*}
\sum_a \big\{ \na \dd\ln\la (U,U) & + &  \dd\ln\la (B_UU) - 2 U{\ln\la }^2  \\
& + & |A^*_{e_a}U|^2 + g((\na_{e_a}B^*)_Ue_a, U) - |B^*_Ue_a|^2 \big\}\,. 
\end{eqnarray*}

Then we have
\begin{eqnarray*}
\sum_a |A^*_{e_a}U|^2  & = & \frac{1}{4}||I||^2 + 2 U(\ln\la )^2\,, \\
g((\na_{e_a}B^*)_Ue_a, U) &  = & g(e_a, \na_{e_a}\mu )\,, \\
\sum_a |B^*_Ue_a|^2 & = & \sum_a g(e_a, \na_UU)^2 = |\mu |^2\,.
\end{eqnarray*}
On noting that $\dd^*\mu^{\flat} = - g(e_a, \na_{e_a}\mu ) + |\mu |^2$,  $\dd^*\Om (U)  =  - \dd^*\mu^{\flat} + |\mu |^2 + \frac{1}{2} |I|^2$ and
$(\Ll_{\mu}g)(U,U) = 2g(\na_U\mu , U) = - 2|\mu |^2$, the formula follows.
\hfill q.e.d.

\bigskip
              
The Ricci tensor evaluated on a horizontal and vertical vector is given by the following expression.
              
\begin{lemma}  \label{lem:mix} 
$$
\Ric (X,U) = X(U(\ln\la )) + \frac{1}{2}\dd^*\Om (X) - U(\ln\la )g(X, \mu ) + \frac{1}{2}\left( \Ll_{\mu}g \right)  (X,U)\,.
$$
\end{lemma}
              
\noindent \emph{Proof}:  By Theorem 11.2.1 of \cite{Ba-Wo}, 
\begin{eqnarray*}  
\Ric (X,U) & = & g(R(X, e_i)e_i, U) = g(R(X, e_a)e_a, U) \\
  & = & g((\na_XA)_{e_a}e_a - (\na_{e_a}A)_Xe_a, U) + g(B^*_Ue_a, I(X, e_a))\,.
\end{eqnarray*}
On calculating, we obtain
\begin{eqnarray*}
 g((\na_XA)_{e_a}e_a, U) - g((\na_{e_a}A)_Xe_a, U) & =  & X(U(\ln\la )) + \frac{1}{2}\dd^*\Om (X)  
- U(\ln\la )g(X, \mu ) \\
& & - \frac{1}{2}g(X, \Ll_{\mu}U) + \frac{1}{2}g(U, \Ll_{\mu}X)\,.
 \end{eqnarray*}
Also $g(B^*_Ue_a, I(X, e_a)) = g(U, [X,\mu ])$ and the formula follows.
 \hfill  q.e.d. 
                 
\bigskip
                 
In order to calculate the Ricci tensor on horizontal vectors, we first establish a useful lemma which we will need later on.
                 
\begin{lemma}  \label{lem:7}  Suppose that, for each $a = 1,2$, the vector field $X_a$ is the horizontal lift of a unit vector field $\ov{X}_a$ on a domain of $N$:  $\dd\phi (X_a) =  \ov{X}_a\circ \phi$.
Then 

{\rm (i)} $\Hh [U, X_a ] = 0$;

{\rm (ii)}  $U\{ g(U, [X_1, X_2])\}  = - \dd\mu^{\flat}(X_1, X_2)$.
\end{lemma}
                 
\noindent \emph{Proof}:  On the one hand $\na\dd\phi (X_a , U) = - \dd\phi (\na_{X_a}U)$.  Whereas
$$
\na \dd\phi (U, X_a) = - \dd\phi (\na_UX_a) + \na_U( \ov{X}_a\circ \phi ) = - \dd\phi (\na_UX_a)\,.
$$
But by the symmetry of the second fundamental form, $\na\dd\phi (X_a, U) = \na\dd\phi (U, X_a)$ and (i) follows.

Now $g(\na_U[X_1, X_2], U) = g([U, [X_1, X_2]]+ \na_{[X_1, X_2]}U, U) = g([U, [X_1, X_2]], U)$.  Furthermore, from the Jacobi identity: $[U, [X_1, X_2]] = [[U, X_1], X_2] + [[X_2, U], X_1]$.  From (i), $\Hh [U, X_a] = 0$; also $\Vv [U, X_a] = - g(X_a, \mu )U$.  Therefore
\begin{eqnarray*} 
[U, [X_1, X_2]] & = & [ - g(X_1, \mu )U, X_2] - [ - g(X_2, \mu )U, X_1] \\
   & =  & - g(X_1, \mu )[U, X_2] + g(X_2, \mu )[U, X_1] 
+ X_2\{ g(X_1, \mu )\} U - X_1\{ g(X_2, \mu )\} U \\
 & = & X_2\{ g(X_1, \mu )\} U - X_1\{ g(X_2, \mu )\} U\,.
\end{eqnarray*}
But $\dd\mu^{\flat} (X_1, X_2) = X_1 g(X_2, \mu ) - X_2 g(X_1, \mu ) - g([X_1,X_2], \mu )$.   The formula follows.

 \hfill q.e.d.   
 
 \begin{lemma} \label{lem:hor}  For horizontal vectors $X, Y$, 
 \begin{eqnarray*}
 \Ric (X,Y) & =  & \Big\{ \la^2 K^N + \Delta \ln\la + \mu (\ln\la ) - \frac{1}{4}|I|^2\Big\} g(X,Y) \\
  & & + \frac{1}{2} \Ll_{\mu}g\,(X,Y) - g(X, \mu )g(Y, \mu )\,.
\end{eqnarray*}
 \end{lemma}
                   
 \noindent \emph{Proof}:  Since both sides are tensorial in $X$ and $Y$, it suffices to verify the formula on a basis.  Let $e_a = \la X_a$, where $X_a$ is the basis taken in Lemma \ref{lem:7}; thus $\{ e_1, e_2\}$ is an orthonormal frame for the horizontal space.  We first of all compute $\Ric (e_1, e_2) = g(R(U, e_1)e_2, U)$.  By Theorem 11.2.1 of \cite{Ba-Wo}, we have
\begin{eqnarray*}
g(R(U, e_1)e_2, U) & =  & g((\na_UA)_{e_1}e_2, U) + g(A^*_{e_1}U, A^*_{e_2}U) \\ 
  & & + g((\na_{e_1}B^*)_Ue_2, U) - g(B^*_Ue_2, B^*_Ue_1) - 2U(\ln\la )g(A_{e_1}e_2, U)\,. 
\end{eqnarray*}
We compute the individual terms in this expression:
$$  
g((\na_UA)_{e_1}e_2, U) = \frac{1}{2} U\{ g([e_1, e_2], U)\}\,.
$$
It is easily checked that $g(A^*_{e_1}U, A^*_{e_2}U)$ vanishes, whereas
$$
g((\na_{e_1}B^*)_Ue_2, U) = g(e_2, \na_{e_1}\mu )\,.
$$
We also have $g(B^*_Ue_2, B^*_Ue_1) = g(e_2, \mu )g(e_1, \mu )$ and $g(A_{e_1}e_2, U) = \frac{1}{2}g([e_1, e_2], U)$.  Combining and applying Lemma \ref{lem:7}, gives the required formula with $X = e_1$ and $Y = e_2$.
                       
On the other hand
\begin{eqnarray*}
\Ric (e_1, e_1) & = & g(R(U, e_1)e_1, U) + ( R(e_2, e_1)e_1, e_2) \\
 & = & K(e_1\wedge U) + K(e_2\wedge e_1)\,, 
\end{eqnarray*}
We now apply Proposition 11.2.2 of \cite{Ba-Wo}, which expresses the sectional curvatures of a semi-conformal submersion, to give the required formula with $X = Y = e_1$.  We omit the details.
\hfill q.e.d. 

\bigskip
                        
\noindent \emph{ Proof of Proposition} \ref{pro:Ricci}:    The formula of the proposition follows on combining Lemmata \ref{lem:vert}, \ref{lem:mix} and \ref{lem:hor}.
\hfill q.e.d. 

\medskip

  We will now equate the expression for the Ricci curvature with the right hand side of (\ref{eq:soliton}).  However, it is useful to decompose the vector field $E$ into its different components.  It is no loss of generality to set $E = - \mu + X + fU$ for some function $f$ and horizontal vector field $X$, in order to cancel the term $\frac{1}{2}\Ll_{\mu}g$ in (\ref{eq:Ricci}).  
\begin{lemma} \label{lem:f}\begin{equation} \label{eq:f}
\Ll_{(fU)}g = - 2f\, U(\ln\la )(g - \theta^2) + 2\dd f\odot \theta + 2 f \, \mu^{\flat}\odot \theta\,.\end{equation}
\end{lemma}

\noindent \emph{Proof}:   It suffices to evaluate $\Ll_{(fU)}g$ on different combinations of vectors, using the fact that $fU$ is tangent to a conformal foliation, so that  $(\Ll_{(fU)}g)(X,Y) = $ $- 2 f U(\ln\la )g(X,Y)$ for all $X,Y\in \Hh$. \hfill q.e.d.

\medskip

We can now write out the fundamental equations for a soliton derived from a semi-conformal map.
\begin{proposition} \label{pro:soliton}  Let $\phi : (M^3, g)\ra (N^2, h)$ be a ($C^{\infty}$) submersive semi-conformal map.  Then $(M^3, g)$ is a Ricci soliton if and only if there is a function $f : M^3 \ra \RR$, a horizontal vector field $X$ and a constant $A$, such that, in the notations defined above,
\begin{eqnarray}  
0 & = & \left\{ \la^2 K^N + \Delta \ln \la + \mu (\ln \la ) - f\, U(\ln \la )\right\} (g - \theta^2) - \frac{1}{4} |I|^2 g   + \frac{1}{2}\Ll_{X} g  + Ag \nonumber \\ 
  &     -  & (\mu^{\flat} + \dd^{\Vv}\ln \la )^2 - (\dd^{\Vv}\ln\la )^2 + \left\{ \dd f + f\, \mu^{\flat} + 2 \dd (U(\ln\la )) + \dd^*\Om\right\}  \odot \theta \,. \label{eq:soliton2}  \end{eqnarray}  
\end{proposition}

\begin{example} \label{ex:wp} {\rm (Warped product solutions)    
Locally, a warped product is of the form $M^3 = N^2 \times J$ with metric $g = (h/\la^2) + \dd t^2$, where $J$ is an open interval in $\RR$ with its coordinate $t$ and  
where $\la = \la (t)$.  We take $\phi$ to be the canonical projection $\phi : M\ra N$.  Then the fibres are geodesic, so that, in addition to being semi-conformal, $\phi$ now has the further property of being a harmonic morphism and our expression (\ref{eq:Ricci}) for the Ricci curvature reduces to a well-known, much simpler case of our formula, see \cite{Ba-Wo}.  Note also that the horizontal distribution is integrable so that $I$ and $\Om$ vanish.  We will suppose that the horizontal component $X$ of the vector field $E$ vanishes and that $K^N$ is constant.  Then the system of equations (\ref{eq:soliton2}) is equivalent to the pair of equations 
\begin{equation} \label{eq:wp} 
\left\{  \begin{array}{lll}  0 & = & \la^2 K^N + \Delta \ln \la - f U(\ln \la ) + A \\  0 & = & U(f) + 2 U(U(\ln\la )) - 2U(\ln\la )^2 + A   
\end{array} \right.  
\end{equation}    

Let $t$ denote a unit speed parameter along the fibres, so that $t = {\rm const.}$ gives the (integrable) horizontal spaces, $U = \pa /\pa t$ and $\theta = \dd t$.  As a consequence of the equations, we also have $f = f(t)$.  We note the special solutions given by $\la \equiv$ const..  Without loss of generality, we may suppose that $\la \equiv 1$.  If now $K^N = 1$, then $A = -1$ and $f^{\prime}= 1$.  This gives the soliton $S^2 \times \RR$ with $E = t (\pa / \pa t)$.  Similarly, if $K^N = -1$, then $A = 1$ and $f^{\prime} = -1$ leading to the soliton $H^2 \times \RR$ with $E = - t(\pa /\pa t)$.

We will now suppose that $\la $ is non-constant and work on a neighbourhood where $\la^{\prime} \neq 0$.  Then, letting $\{ e_a\}$ denote a local orthonormal frame, by Lemma \ref{lem:Hmc}, $\Delta \ln \la  =  \Tr \na \dd\ln\la = \dd\ln\la (\na_{e_a}e_a) + U(U(\ln\la ))$ $=  - 2\{ (\ln\la )^{\prime}\}^2 + (\ln\la )^{\prime\prime}$.   

The system (\ref{eq:wp}) now becomes the pair of ordinary differential equations:    
\begin{equation} \label{eq:wp2}    
\left\{  \begin{array}{lll}     0  & = & \la^2 K^N + (\ln\la )^{\prime\prime} - 2 \{ (\ln\la )^{\prime}\}^2 - f (\ln\la )^{\prime}  + A \\     0 & = & f^{\prime} + 2 (\ln\la )^{\prime\prime} - 2 \{ (\ln\la )^{\prime}\}^2  + A      
\end{array}\right.     
\end{equation}    
 We can solve these to express $f$ in terms of $\la$:    
 \begin{equation} \label{eq:wp-f}     f = \frac{\la^{\prime\prime} + A \la + K^N \la^3}{\la^{\prime}} - \frac{3\la^{\prime}}{\la}\,.     
\end{equation}     
Substituting back, we obtain the following 3rd order ordinary differential equation:     
\begin{equation} \label{eq:wp-3rdorder}\la^{\prime\prime\prime} \la^{\prime}\la^2 - \la\la^{\prime\prime}(\la^{\prime\,2} + \la\la^{\prime\prime})+ A\la^2 (2\la^{\prime\,2} - \la\la^{\prime\prime}) + K^N\la^4(3\la^{\prime\,2}- \la\la^{\prime\prime}) - \la^{\prime\,4} = 0\,.
\end{equation}
Thus (\ref{eq:wp-f}) and (\ref{eq:wp-3rdorder}) are equivalent to (\ref{eq:wp2}). 
It is easily checked that the solution has constant curvature if and only if 
$(\ln\la )^{\prime\prime} - \la^2K^N = 0$. 

If we look for solutions of the form 
\begin{equation} \label{eq:elliptic}(\la^{\prime})^2 = 2F(\la )
\end{equation} 
for some function $F$,  then (\ref{eq:wp-3rdorder}) becomes
$$
2\la^2 FF^{\prime\prime} - \la F^{\prime}(2F + \la F^{\prime}) + A \la^2 (4F - \la F^{\prime}) + K^N\la^4 (6F - \la F^{\prime}) - 4F^2 = 0\,.
$$
It is easy to see that the only polynomial solutions are given by 
$$
F = \frac{K^N \la^4}{2} + \frac{A\la^2}{4}\,.
$$
Then (\ref{eq:elliptic}) can be integrated explicitly, however all these solutions have constant curvature.
On the other hand, the exceptional solutions when $A = K^N = 0$, given by $F(\la ) = B\la^k$ with $k = 2(1\pm \sqrt{2})$ are of non-constant curvature.  Choosing $B = 2$, this gives $\la (t) = t^{\pm 1/\sqrt{2}}$.  We may express the corresponding (incomplete) metric in the form $g = t^{\pm\sqrt{2}}(\dd x^2 + \dd y^2) + \dd t^2$.  }
\end{example} 

\section{Prescribing the mean-curvature of the fibres} \label{sec:mc-fibres}

We wish to characterize those semi-conformal maps $\phi : M^3\ra N^2$ which can be recovered from data on the surface $N^2$.  Consider the $2$-form $\Om = \dd \theta$.  Then by Lemma \ref{lem:app-1}, 
$$
\Om = \wt{\Om} - \mu^{\flat}\wedge \theta\,,
$$
where $\wt{\Om} = \Om \circ \Hh = - g(U, I)$, with $I$ the integrability tensor $I(X,Y) = \Vv [X,Y]$.  Say that a $1$-form $\om$ is \emph{closed relative to a distribution} $D$ if and only if $\dd\om (X,Y) = 0$ for all $X,Y\in D$.

\begin{lemma} \label{lem:rel-closed}
The $2$-form $\wt{\Om}$ is basic if and only if $\mu^{\flat}$ is closed relative to the horizontal distribution, i.e. $\dd \mu^{\flat}(X,Y) = 0$ for all $X,Y\in \Hh$.  In particular, this is the case whenever $\mu$ is the gradient of a function.
\end{lemma}

\noindent \emph{Proof}:  If $X,Y$ are basic horizontal vector fields, then, as in Lemma \ref{lem:7}(i), $\Hh\Ll_UX = \Hh\Ll_UY=0$ and
\begin{eqnarray*}
(\Ll_U\wt{\Om})(X,Y) & = & U\left( \wt{\Om}(X,Y)\right) - \wt{\Om}(\Ll_UX, Y) - \wt{\Om}(X, \Ll_UY) \\
 & = & U(\wt{\Om}(X,Y)) = U(\Om (X,Y))\,.
\end{eqnarray*}
It follows that $\wt{\Om}$ is basic if and only if $U\{ \Om (X,Y)\} = 0$ for all \emph{basic} horizontal vector fields $X,Y$.  But by Lemma \ref{lem:7}(ii), $U\{ \wt{\Om}(X,Y)\} = \dd \mu^{\flat}(X,Y)$.  The result follows.  \hfill  q.e.d.

\begin{corollary} If $\mu = \Hh \grad \ln \rho$ for some function $\rho$.  Then $\mu^{\flat}$ is closed relative to the horizontal distribution if and only if {\rm either} $\wt{\Om} \equiv 0$ (the horizontal distribution is integrable), {\rm or} $U(\ln\rho ) = 0$, i.e. $\rho$ is a basic function.
\end{corollary}

\noindent \emph{Proof}:  Let $X,Y$ be horizontal vector fields.  Then
\begin{eqnarray*}
\dd\mu^{\flat}(X,Y) & = & X(Y(\ln\rho )) - Y(X(\ln\rho )) - g(\Hh \grad \ln \rho , [X,Y]) \\
 & = & g(\Vv \grad \ln \rho , [X,Y]) = - U(\ln \rho ) \wt{\Om} (X,Y)\,.
\end{eqnarray*}
The result follows.  \hfill q.e.d. 

\medskip

Suppose that $\phi : (M^3, g) \ra (N^2, h)$ is a semi-conformal map with dilation $\la$ and with associated integrability $2$-form $\wt{\Om }$ basic.  We will suppose further that $M^3$ is a  domain whose $2$-cohomology vanishes and on which $\phi$ is submersive with connected fibres.  This is no loss of generality, since otherwise we work locally on an open set with these properties.  Then we may write $\wt{\Om} = \phi^*\ov{\Om}$ for some $2$-form $\ov{\Om}$ on $N$.  Clearly $\wt{\Om}$ is closed, since $\dd\wt{\Om} = d\phi^*\ov{\Om} = \phi^*\dd\ov{\Om}$.  Since $H^2(M, \RR ) = 0$, we have $\wt{\Om} = \dd\wt{\theta}$ for some $1$-form $\wt{\theta}$ (defined up to addition of a derivative).  Set 
$$
\wt{g} = \frac{\phi^*h}{\la^2} + \wt{\theta}^2\,.
$$

\begin{proposition} \label{prop:semi-har}
The map $\phi : (M^3, \wt{g}) \ra (N^2, h)$ is both semi-conformal and harmonic, equivalently, $\phi$ is a harmonic morphism with respect to the metric $\wt{g}$.
\end{proposition}

\noindent \emph{Proof}:  For a semi-conformal map onto a surface, harmonicity is equivalent to the fibres being minimal (\cite{Ba}).  Let $\wt{U}$ be a unit vertical vector field with respect to $\wt{g}$.  Then it is easily checked that the fibres are geodesic if and only if $\Ll_{\wt{U}}\wt{\theta} = 0$ (see \cite{Ba-Wo}).  But 
$$
\Ll_{\wt{U}}\wt{\theta} = \dd (\wt{\theta}\rfloor \wt{U}) + \dd\wt{\theta}\rfloor \wt{U} = \wt{\Om}\rfloor \wt{U} = 0\,.
$$
\hfill  q.e.d. 

\medskip

We therefore see how, to a semi-conformal map $\phi : (M^3, g)\ra (N^2, h)$ with $\mu^{\flat}$ closed relative to the horizontal distribution, we can associate a harmonic morphism $\phi : (M^3, \wt{g})\ra (N^2, h)$.  We now consider the converse problem:  given a harmonic morphism $\phi : (M^3, \wt{g}) \ra (N^2, h)$, construct a semi-conformal map $\phi : (M^3, g)\ra (N^2, h)$ with fibres having \emph{prescribed} mean curvature $\mu$. 

\begin{proposition} \label{prop:prescribed}
Let $\phi : (M^3, \wt{g}) \ra (N^2, h)$ be a submersive harmonic morphism with dilation $\la$ and integrability $2$-form $\wt{\Om}$.  Let $\eta$ be a given $1$-form and let $\theta$ be a  $1$-form satisfying the exterior differential equation
\begin{equation} \label{eq:ext-bis}
\dd\theta +  \eta \wedge \theta = \wt{\Om}\,.
\end{equation}
Suppose that $\theta (\wt{U})\neq 0$ at every point, where $\wt{U}$ is a unit vertical vector with respect to $\wt{g}$.  Let $g = \frac{\phi^*h}{\la^2} + \theta^2$.  Then $\phi : (M^3, g) \ra (N^2, h)$ is semi-conformal with fibres having mean curvature $\Hh(\eta^{\sharp})$. Furthermore, if $\eta$ is the derivative of a basic function, then we can always solve {\rm (\ref{eq:ext-bis})} on any domain $M^3$ with $H^2(M^3, \RR)=0$.
\end{proposition}

\noindent \emph{Proof}:  Clearly $\phi : (M^3, g) \ra (N^2, h)$ is semi-conformal.  Set $\Om = \dd\theta$.  Then by Lemma \ref{lem:app-1}, $\Om (X,U) = - \mu^{\flat}(X) = - g(\mu , X)$, where $X$ is horizontal (with respect to $g$) and $\mu$ denotes the mean-curvature of the fibres.  But since $U$ and $\wt{U}$ are colinear (both being in the kernel of $\dd\phi$), $\wt{\Om} \rfloor U = 0$ and by (\ref{eq:ext-bis}), $\Om (X,U) = - \eta (X)$, i.e. $\eta (X) = g(\mu , X)$ as required.  

For the last part of the proposition, write $\eta = \dd \ln \rho$, where $\rho = \ov{\rho}\circ \phi$.  Then equation (\ref{eq:ext-bis}) has the form
$$
\dd (\rho \theta ) = \rho \wt{\Om}\,.
$$
Then we can solve the equation $\dd\om = \rho \wt{\Om}$ on a domain with vanishing $2$-cohomology if and only if $\dd (\rho \wt{\Om}) = 0$.  But $\rho \wt{\Om} = \phi^*(\ov{\rho}\ov{\Om})$, which is clearly closed. 
\hfill q.e.d.

\medskip 

In order to construct a semi-conformal map  from data on $N^2$ we proceed as follows.  Let $(N^2, h)$ be a Riemannian surface with $H^2(N^2, \RR ) = 0$.  Let $\ov{\la}, \ov{\rho} : N^2\ra \RR$ be smooth positive functions and let $\ov{\Om}$ be a smooth $2$-form on $N^2$ (equivalently, we may let $\ov{\si} : N^2 \ra \RR$ be a smooth function and write $\ov{\Om} = \ov{\si} \mu^N$, where $\mu^N$ is the volume form on $N^2$).  Set $M^3 = N^2 \times (- \delta , \delta )$ for some $\delta >0$ and let $\phi : M^3 \ra N^2$ be the projection.  Define $\wt{\Om} = \phi^*\ov{\Om}$.  Then $\dd\wt{\Om} = 0$ so that by vanishing $2$-cohomology, we have $\wt{\Om} = \dd\wt{\theta}$ for some $1$-form $\wt{\theta}$.  On writing $t$ for the coordinate of $(- \delta , \delta )$, by replacing $\wt{\theta}$ with $\wt{\theta} + R\dd t$ if necessary ($R$ a sufficiently large constant) we may suppose that $\wt{\theta}(\pa /\pa t) \neq 0$ at every point.  Set $\la = \ov{\la}\circ \phi$, $\rho = \ov{\rho}\circ \phi$ and let $\wt{g} = \frac{\phi^*h}{\la^2} + \wt{\theta}^2$.
Then $\phi : (M^3, \wt{g}) \ra (N^2, h)$ is a harmonic morphism.  In fact, since the gradient of the dilation is horizontal, it is a harmonic morphism of \emph{Killing type}, that is the fibres are the integral curves of a Killing vector field (see \cite{Br-1}).  Now let $\theta$ solve the exterior differential equation
$$
\dd\theta + \dd\ln\rho \wedge \theta = \wt{\Om}\,.
$$
By Proposition \ref{prop:prescribed}, if we set $g = \frac{\phi^*h}{\la^2} + \theta^2$, then provided $\theta (\pa /\pa t) \neq 0$, the metric $g$ is positive definite and $\phi : (M^3, g) \ra (N^2, h)$ is a semi-conformal map with fibres having mean curvature $\grad \ln \rho$.

In the next section, we will establish conditions on the functions $\ov{\la}, \ov{\rho}, \ov{\si}$ defined above, which are equivalent to the property that $(M^3, g)$ be a Ricci soliton.

\section{Solitons constructed from data on a surface} \label{sec:killing}

Our aim in this section is to establish Theorem \ref{th:main}.  The ansatz is derived from the construction of a semi-conformal mapping whose dilation and  mean curvature of its fibres are both basic.   

Let $\phi : M^3\ra N^2$ be a semi-conformal submersion with connected fibres, basic dilation $\la = \ov{\la}\circ \phi$ and whose fibres have mean curvature which is the gradient of a basic function: $\mu = \grad \ln \rho$ ($\rho = \ov{\rho}\circ \phi$).
A special type of horizontal vector field $X$ is given by the gradient of a basic function, that is $X = \grad \ln\nu$, where $\nu = \ov{\nu}\circ \phi$ with $\ov{\nu}:N\ra \RR$.  We will now make this assumption on $X$; in particular this implies that $\Ll_Xg = 2\na \dd \ln\nu$.  This need not imply that a corresponding soliton is of \emph{gradient type}; for this we require that $E$ be the gradient of a function, which depends essentially on the function $f$.  Since $\dd \phi (\grad \ln \nu )= \ov{\la}^2\grad \ln \ov{\nu}\circ \phi$, we have the correspondence $\ov{X} = \ov{\la}^2\grad \ln \ov{\nu}$.  Thus with reference to Remark 2 of the Introduction, we need to make the substitution $\ov{Y} = \ov{X}/\ov{\la}^2$.  

 Equation (\ref{eq:soliton2}) for a soliton now has the form:
\begin{equation} \label{eq:Ksoliton}
\{ \la^2 K^N + \Delta^M \ln\la + \mu (\ln\la ) \} g|_{\Hh\times \Hh} - \frac{1}{4} ||I||^2g + Ag + \na \dd \ln\nu - (\mu^{\flat})^2 + \{ \dd f + f\,u^{\flat} + \dd^*\Om\}\odot \theta = 0\,.
\end{equation}
An immediate necessary condition for this to be satisfied is the requirement that $\na \dd \ln\nu - (\mu^{\flat})^2$ be umbilic on $\Hh$, i.e. $\na \dd \ln\nu |_{\Hh\times \Hh} - (\mu^{\flat})^2= \al \, g|_{\Hh\times \Hh}$ for some function $\al : M\ra \RR$.  We will first of all investigate this condition more fully.  For a covariant tensor $S$, we will write $S\rfloor X$ for its contraction with a vector $X$, so that $(S\rfloor X)(Y_1, \ldots , Y_k) = S(X,Y_1, \ldots , Y_k)$.

\begin{lemma} \label{lem:umbilicF}  Let $F = \ov{F}\circ \phi : M\ra \RR$ be any smooth basic function, then
\begin{equation} \label{eq:pullback}
\na\dd F = \phi^*\left\{ \na \dd \ov{F} + 2\dd\ln\ov{\la}\odot \dd\ov{F} - h(\grad \ln\ov{\la}, \grad \ov{F})h\right\} + \left( \Om \rfloor \grad F\right) \odot\theta\,,
\end{equation}
with $\Om \rfloor \grad F = \phi^*(\ov{\la}^2\ov{\Om}\rfloor \grad \ov{F})$.
\end{lemma}

\noindent \emph{Proof}:  First let $X,Y$ be horizontal vectors and set $\ov{X} = \dd\phi (X)$ etc..  Let $\{ \ov{X}_a\}_{a=1,2}$ be an orthonormal frame on a domain of $N$ and write $X_a$ for the horizontal lift of $\ov{X}_a$.  Note that $\{ \la X_a\}_{a=1,2}$ is an orthonormal frame for $\Hh$.  Then
$$
\na\dd F(X,Y) = - \dd F(\na^M_XY) + X(Y(F)) = - \dd\ov{F}(\dd\phi (\na_X^NY)) + \ov{X}(\ov{Y} (\ov{F}))\circ\phi\,.
$$
But 
\begin{eqnarray}
\dd\phi (\na^M_XY) & = & \dd\phi (g(\na^M_XY, \la X_a)\la X_a) \nonumber \\
 & = & \frac{\la^2}{2}\Big\{ Xg(Y, X_a)+Yg(X, X_a) - X_ag(X,Y) \nonumber \\
  &  &  \quad - g(X,[Y, X_a]) - g(Y, [X, X_a]) + g(X_a,[X,Y]) \Big\} \ov{X}_a \nonumber \\
 & = & \frac{\la^2}{2}\Big\{ X\left( \frac{1}{\la^2}h(\ov{Y}, \ov{X}_a)\right) +Y\left( \frac{1}{\la^2}h(\ov{X}, \ov{X}_a)\right) - X_a\left( \frac{1}{\la^2}h(\ov{X},\ov{Y})\right) \nonumber \\
  &  &  \quad - \frac{1}{\la^2}h(\ov{X},[\ov{Y}, \ov{X}_a]) - \frac{1}{\la^2}h(\ov{Y}, [\ov{X}, \ov{X}_a]) + \frac{1}{\la^2}h(\ov{X}_a,[\ov{X},\ov{Y}]) \Big\} \ov{X}_a \nonumber \\
  & = & \na^N_{\ov{X}}\ov{Y} - X(\ln\la )h(\ov{Y}, \ov{X}_a)\ov{X}_a \nonumber \\
 & & \qquad - Y(\ln\la )h(\ov{X}, \ov{X}_a)\ov{X}_a + X_a(\ln\la )h(\ov{X}, \ov{Y})\ov{X}_a. \label{eq:proj-cov}
\end{eqnarray}
The formula evaluated on horizontal vectors follows.

Note that $\na \dd F (U,U) = - \dd F(\na_UU) = - \dd F(\mu ) = (\Om \rfloor \grad F)(U)$.  On the other hand 
\begin{eqnarray*}
\na \dd F(X_a, U) & = & - \dd\ov{F}(\dd\phi (\na_{X_a}U) = - \la^2\dd\ov{F} \left( \dd\phi (<\na_{X_a}U, X_b>X_b)\right) \\
  & = & - \la^2<\na_{X_a}U, X_b>\ov{X}_b(\ov{F}) \\
 & = & - \frac{\la^2}{2} \Big\{ X_ag(U, X_b) + Ug(X_a, X_b) - X_b g(X_a, U) \\
 & & \quad - g(X_a, [U, X_b]) - g(U, [X_a, X_b]) + g(X_b, [X_a, U])\Big\} \ov{X}_b(\ov{F})\,.
\end{eqnarray*}
Now $\Hh [U, X_a] = 0$ and $g(X_a, X_b) = \delta_{ab}/\la^2$, but since $U(\la )=0$, then $Ug(X_a, X_b)=0$, hence the above expression equals
$$
\frac{\la^2}{2}g(U, [X_a, X_b])\ov{X}_b(\ov{F}) = - \frac{\la^2}{2}\Om (X_a, X_b)\ov{X}_b(\ov{F}) = \frac{\la^2}{2}\ov{\Om}(\ov{X}_b(\ov{F})\ov{X}_b, \ov{X}_a).
$$
Finally, it is easily checked that $\dd\phi (\grad F) = \la^2\grad \ov{F} \circ \phi$.
\hfill q.e.d.  

\medskip

The following corollary in fact holds for any semi-conformal  map $\phi : M^m\ra N^2$ without further hypotheses, however, it is an easy consequence of (\ref{eq:pullback}).

\begin{corollary} \label{cor:lap}
Let $F = \ov{F}\circ \phi : M^3 \ra \RR$ be any smooth basic function, then
$$
\Delta^MF+\mu (F) = \la^2\Delta^N\ov{F}\circ \phi\,.
$$
\end{corollary}

\noindent \emph{Proof}:  We have $\Delta^MF = \Tr \na \dd F = \Tr_{\Hh}\na\dd F + \na \dd F (U,U)$.  But $\na \dd F(U,U)$ $ = - \dd F(\na_UU)$ $ = - \mu (F)$.  On taking the trace over the horizontal space in (\ref{eq:pullback}), the formula follows.  \hfill q.e.d.

\medskip

We can now characterize the umbilicity condition discussed above.

\begin{corollary} \label{cor:umbilicF}  
$\na \dd \ln\nu |_{\Hh\times \Hh} - (\mu^{\flat})^2$ is umbilic on $\Hh$ if and only if
\begin{equation} \label{eq:umbilic2}
\na \dd \ln\ov{\nu} + 2 \dd\ln\ov{\la} \odot \dd\ln \ov{\nu} - (\dd \ln \ov{\rho})^2 = \al h\,,
\end{equation}
for some function $\al : N \ra \RR$.
\end{corollary}

\noindent \emph{Proof}:  For the mean curvature, we have $\mu^{\flat} = \dd\ln\rho = \phi^*(\dd\ln\ov{\rho})$, so that from (\ref{eq:pullback}),
$$
\na \dd \ln\nu |_{\Hh\times \Hh} - (\mu^{\flat})^2 = \phi^*\{ \na \dd\ln \ov{\nu} + 2\dd \ln\ov{\la} \odot \dd\ln \ov{\nu} - h(\grad \ln\ov{\la}, \grad\ln \ov{\nu})h - (\dd\ln\ov{\rho})^2\}\,.
$$
We require the bracket to be proportional to $h$, which is the assertion of the corollary.

\hfill q.e.d. 

\medskip

By taking traces, we have the following consequence.

\begin{corollary} \label{cor:umbilic-trace}
If $\na \dd \ln\nu|_{\Hh\times \Hh} - (\mu^{\flat})^2$ is umbilic on $\Hh\times \Hh$, then
\begin{eqnarray*}
\na \dd \ln\nu |_{\Hh\times \Hh} - (\mu^{\flat})^2  & = & \frac{1}{2}\left( \Delta^M \ln \nu + \mu (\ln\nu ) - |\mu |^2 \right) g|_{\Hh\times \Hh} \\
 & = & \frac{1}{2}\phi^*\left\{ (\Delta^N\ln\ov{\nu} - |\grad\ln\ov{\rho}|^2)h\right\}\,.
\end{eqnarray*}
\end{corollary}

\begin{lemma} \label{lem:I}
Let $\wt{\Om}$ be the integrability $2$-form associated to $\phi$ and let $\wt{g}$ be the corresponding metric (cf. Section {\rm \ref{sec:mc-fibres}}).  Then $||\wt{\Om}||^2_{\wt{g}} = ||I||_g^2$.  In particular, $||I||^2_g$ is a basic function.
\end{lemma} 

\noindent \emph{Proof}:  Let $Y$ be horizontal and basic with respect to $g$: $\dd\phi (Y) = \ov{Y}, g(Y,U)=0$.  Then, recalling that $\wt{g} = \frac{\phi^*h}{\la^2} + \wt{\theta}^2$, it follows that $\wt{Y} = Y - \frac{\wt{\theta}(Y)}{\wt{\theta}(U)}U$ is horizontal and basic with respect to $\wt{g}$.  With the same notations, let $\{ \ov{X}_a\}$ be an orthonormal basis on $N$, so that $\{ \la X_a\}$ is an orthonormal basis for $\Hh$ on $(M, g)$ and $\{ \la \wt{X}_a\}$ is an orthonormal basis for $\wt{\Hh}$ on $(M,\wt{g})$.  But $\wt{\Om}\rfloor U = 0$, so that
\begin{eqnarray*}
||\wt{\Om}||^2_{\wt{g}} & = & 2 \wt{\Om}(\la \wt{X}_1, \la \wt{X}_2)^2 = 2\la^4\wt{\Om}(\wt{X}_1, \wt{X}_2)^2 = 2\la^4\wt{\Om}(X_1, X_2)^2 \\
 & = & 2 \la^4\Om (X_1, X_2)^2 = 2 \Om (\la X_1, \la X_2)^2 = ||I||^2_g.
\end{eqnarray*}
For the last part:
$$
||\wt{\Om}||^2_{\wt{g}} = 2\la^4\wt{\Om}(\wt{X}_1, \wt{X}_2)^2 = 2\la^4 \phi^*\ov{\Om}(\wt{X}_1, \wt{X}_2)^2 = 2\la^4\ov{\Om}(\ov{X}_1, \ov{X}_2)^2\circ\phi = \la^4||\ov{\Om}||_h^2\,.
$$
But since $\la = \ov{\la} \circ \phi$ is basic, so is $||I||_g^2$.  
\hfill q.e.d. 

\medskip 

It will be convenient in what follows to write $\frac{1}{4}||I||^2 = \psi = \ov{\psi} \circ \phi$. Let us review the soliton equation (\ref{eq:Ksoliton}) under the assumption that the umbilicity condition (\ref{eq:main})(u) is satisfied.  By the above calculations it takes the form:
\begin{eqnarray} 
0 & = & \left\{ \ov{\la}^2K^N + \ov{\la}^2 \Delta^N\ln\ov{\la} + \frac{\ov{\la}^2}{2}\left( \Delta^N\ln\ov{\nu} - |\grad \ln\ov{\rho}|^2\right) - \ov{\psi} + A \right\} g_{\Hh\times \Hh} \nonumber \\
 &  &  \qquad  + \left\{ \dd f + f\mu^{\flat} + \Om \rfloor \grad\ln\nu + \dd^*\Om - \psi \theta + A\theta \right\} \odot \theta\,, \label{eq:soliton-neat}
\end{eqnarray}
with its component parts:
\begin{equation} \label{eq:soliton-parts}
\left\{ \begin{array}{ll}
{\rm (i)} & K^N + \Delta^N \ln\ov{\la} + \frac{1}{2}\left( \Delta^N\ln\ov{\nu} - |\grad \ln\ov{\rho}|^2\right) - \frac{\ov{\psi} - A}{\ov{\la}^2} = 0  \\
{\rm (ii)} & \dd f + f\mu^{\flat} + \Om \rfloor \grad\ln\nu + \dd^*\Om - (\psi - A)\theta = 0 \,.
\end{array} \right.
\end{equation}
We will now study (\ref{eq:soliton-parts})(ii).  Recall that $\Om = \wt{\Om} - \mu^{\flat}\wedge \theta = \phi^*\ov{\Om} - \dd\ln\rho \wedge \theta$.  The following lemma is easily established.

\begin{lemma} \label{lem:wedge}
 For two $1$-forms $\om , \eta$, we have
$$
\dd^*(\om\wedge\eta ) = \dd^*\om \cdot \eta - \dd^*\eta \cdot \om - [\om^{\sharp}, \eta^{\sharp}]^{\flat}\,.
$$
In particular
$$
\dd^*(\mu^{\flat}\wedge \theta ) = (\dd^*\mu^{\flat} - |\mu |^2)\theta = - (\Delta^M\rho - |\mu |^2)\theta = - \la^2(\Delta^N\ln\ov{\rho})\circ\phi \cdot \theta\,.
$$
\end{lemma}

\begin{lemma} \label{lem:dstar-omega}
$$
\dd^*\wt{\Om} = \la^2\phi^*\left\{ \dd^*\ov{\Om} +\ov{\Om}\rfloor \grad \ln (\ov{\rho}\ov{\la}^{-2})\right\} + 2\psi \theta\,.
$$
\end{lemma}

\noindent \emph{Proof}:   
Let $\{\ov{X}_a\}$ be an orthonormal frame on a domain of $N$ and let $X_a$ be the horizontal lift of $\ov{X}_a\ (a = 1,2)$.  Set $e_a = \la X_a$ and let $e_3 = U$, so that $\{ e_1,e_2,e_3\}$ is an orthonormal frame on a domain of $M$.  Then
\begin{eqnarray*}
 \dd^*\wt{\Om} (X_b) & = & - (\na_{e_i}\wt{\Om} )(e_i, X_b) 
  =  - \la^2(\na_{X_a}\wt{\Om})(X_a, X_b) - (\na_U\wt{\Om})(U, X_b) \\
  & = & - \la^2X_a\{ \wt{\Om} (X_a, X_b)\}  + \la^2\wt{\Om} (\na_{X_a}X_a, X_b) + \la^2\wt{\Om} (X_a, \na_{X_a}X_b)+ \wt{\Om}(\na_UU, X_b) \\
   & = & - \la^2X_a\{ \ov{\Om} (\ov{X}_a, \ov{X}_b)\circ\phi\}  + \la^2\ov{\Om} (\dd\phi(\na_{X_a}X_a), \ov{X}_b) \\
 &  &  \qquad \qquad \qquad + \la^2\ov{\Om} (\ov{X}_a, \dd\phi(\na_{X_a}X_b)) + \ov{\Om}(\dd\phi (\mu ), \ov{X}_b)\,.
   \end{eqnarray*}
We now employ the expression (\ref{eq:proj-cov}) for $\dd\phi (\na_{X_a}X_b)$, to give
\begin{equation} \label{eq:d*hor}
 \dd^*\wt{\Om} (X_b) = \la^2 \left\{ \dd^*\ov{\Om}(\ov{X}_b) - 2\ov{\Om}(\grad \ln\ov{\la}, \ov{X}_b) + \ov{\Om}(\grad \ln \ov{\rho}, \ov{X}_b)\right\}\circ \phi\,.
\end{equation}

On the other hand
\begin{eqnarray*}
\dd^*\wt{\Om}(U) & = & - (\na_{e_i}\wt{\Om})(e_i, U) = - \la^2(\na_{X_a}\wt{\Om})(X_a, U) \\
 & = & \la^2\wt{\Om}(X_a, \na_{X_a}U) = \la^2 \wt{\Om}(X_a, g(\na_{X_a}U, e_b)e_b) \\
 & = & - g(U, \na_{e_a}e_b)\wt{\Om}(e_a, e_b) = - \frac{1}{2}g(U, [e_a, e_b])\wt{\Om}(e_a, e_b) 
 =  \frac{1}{2}||I||^2\,.
\end{eqnarray*}
\hfill q.e.d. 

\medskip

We can now express equation (\ref{eq:soliton-parts})(ii) in the form:
\begin{eqnarray} 
\dd (f\rho ) & = & -\rho \la^2 \phi^*\left\{ \dd^*\ov{\Om}+ \ov{\Om}\rfloor \grad \ln (\ov{\rho}\ov{\nu}\ov{\la}^{-2}) \right\} \nonumber \\
    &  &  \quad - \rho\la^2\left\{ \Delta^N\ln\ov{\rho} - h(\grad \ln\ov{\rho}, \grad\ln \ov{\nu}) + \frac{A+\ov{\psi}}{\ov{\la^2}}\right\}\theta \,. \label{eq:parts2-bis}
\end{eqnarray}
Written in this way, we see that, on a simply connected domain, we can find a function $f$ such that (\ref{eq:soliton-parts})(ii) is satisfied if and only if the derivative of the right-hand side vanishes, that is, if and only if
\begin{eqnarray*}
 & \dd \ln (\la^2\rho )\wedge \phi^*\left\{ \dd^*\ov{\Om} +\ov{\Om}\rfloor \grad \ln (\ov{\rho}\ov{\nu}\ov{\la}^{-2}) \right\} \\
 & + \left( \Delta^N\ln\ov{\rho} - h(\grad \ln\ov{\rho}, \grad\ln \ov{\nu}) + \frac{A+\ov{\psi}}{\ov{\la}^2}\right) \circ \phi \cdot \dd\ln\ov{\la}^2 \wedge \theta \\
 & + \phi^*\Delta \ov{\Om} + \phi^* \dd \left\{ \ov{\Om}\rfloor \grad \ln (\ov{\rho}\ov{\nu}\ov{\la}^{-2}) \right\} \\
  & + \dd\left\{ \Delta^N\ln\ov{\rho} - h(\grad \ln\ov{\rho}, \grad \ln\ov{\nu}) + \frac{A+\ov{\psi}}{\ov{\la}^2}\right\} \circ\phi \wedge \theta \\
  & + \left\{ \Delta^N\ln\ov{\rho} - h(\grad \ln\ov{\rho}, \grad \ln\ov{\nu}) + \frac{A+\ov{\psi}}{\ov{\la}^2}\right\}\circ \phi\cdot \phi^*\ov{\Om} = 0.
\end{eqnarray*}
Here, $\Delta \ov{\Om}$ is the Laplacian on forms: $\Delta = \dd\dd^*+\dd^*\dd$.
On simplifying and separating into components, we have equivalence with the following pair of equations on $N$\,:
\begin{eqnarray*}
 {\rm (a)}  & & \dd \ln (\ov{\la}^2\ov{\rho}) \wedge \left\{ \dd^*\ov{\Om} + \ov{\Om}\rfloor \grad \ln (\ov{\rho}\ov{\nu}\ov{\la}^{-2}) \right\} 
 + \dd\left\{ \ov{\Om}\rfloor \grad \ln (\ov{\rho}\ov{\nu}\ov{\la}^{-2}) \right\} + \Delta \ov{\Om} \\
  & & \quad + \left\{ \Delta^N\ln\ov{\rho} - h(\grad \ln\ov{\rho}, \grad\ln \ov{\nu}) + \frac{A+\ov{\psi}}{\ov{\la}^2} \right\} \ov{\Om} = 0  \\
{\rm (b)} & & \left\{\Delta^N\ln\ov{\rho} - h(\grad \ln\ov{\rho}, \grad\ln \ov{\nu}) + \frac{A+\ov{\psi}}{\ov{\la}^2}\right\} \dd\ln\ov{\la}^2  \\
  & & \quad +\dd \left\{ \Delta^N\ln\ov{\rho} - h(\grad \ln\ov{\rho}, \grad \ln\ov{\nu}) + \frac{A+\ov{\psi}}{\ov{\la}^2}\right\} = 0\,.
\end{eqnarray*}

In order to simplify these, it is useful to write $\ov{\Om} = \ov{\si} \mu^N$, where $\mu^N$ is the volume form on $(N,h)$.  Then it is easily checked that the following relation holds:
\begin{equation} \label{eq:rel}
2\ov{\psi} = \ov{\la}^4\ov{\si}^2\,.
\end{equation}
The following identities are useful.
\begin{lemma} \label{lem:ids}
Let $\al, \be$ be arbitrary smooth functions on $N$, then

{\rm (i)}  $\dd \al\wedge (\mu^N\rfloor \grad \be ) = h(\grad \al , \grad \be )\mu^N\,$;

{\rm (ii)}  $ \dd (\mu^N\rfloor \grad \al ) = \Delta^N\al \cdot \mu^N$\,;

{\rm (iii)}  $\dd^*(\al \mu^N) = - \mu^N\rfloor \grad \al$\,.
\end{lemma}

On applying these to equation (a) above, we obtain
\begin{eqnarray*}
\Big\{ \Delta^N\ln (\ov{\rho}^2\ov{\nu}\ov{\la}^{-2}\ov{\si}^{-1})  & + & |\grad\ln\ov{\rho}|^2 - |\grad\ln (\ov{\la}^2\ov{\si})|^2 \\
 &  & + h(\grad \ln (\ov{\la}^2\ov{\si}), \grad \ln\ov{\nu}) 
 +   \frac{A+\ov{\psi}}{\ov{\la}^2}\Big\} \ov{\Om} = 0\,.
\end{eqnarray*}
From (\ref{eq:rel}) we can eliminate $\ov{\si}$.  Then, \emph{either} $\ov{\Om} = 0$, i.e. $\ov{\psi}=0$ and the horizontal distribution is integrable, \emph{or} equation (ii)(a) of Theorem \ref{th:main} holds.  Equation (ii)(b) of the theorem is immediate from (b) above.  This establishes Theorem \ref{th:main}.

\medskip

We now isolate the condition of constant curvature.

\begin{proposition} \label{prop:cc}
Let $(M^3, g)$ be a Ricci soliton as in Theorem {\rm \ref{th:main}}.  Then $(M^3, g)$ has constant curvature if and only if both

{\rm (i)}  {\rm either} $\ov{\psi} = 0$ ($\Hh$ integrable) {\rm or} $\ov{\rho}^2/\ov{\psi}^{-1/2}=$ const.;  and  
\begin{eqnarray*}
{\rm (ii)}\ 
\Big\{ K^N + \Delta^N\ln (\ov{\la}\ov{\rho}^{-1})  & + & |\grad \ln\ov{\rho}|^2 - h(\grad\ln\ov{\la}, \grad\ln\ov{\rho}) - \frac{2\ov{\psi}}{\ov{\la}^2}\Big\} h\\
& & + \na\dd\ln\ov{\rho} + 2\dd\ln\ov{\la}\odot\dd\ln\ov{\rho} - (\dd\ln\ov{\rho})^2   = 0 \,.
\end{eqnarray*}
are satisfied.
\end{proposition}

\noindent \emph{Proof}:  The manifold $M^3$ has constant curvature if and only if, for all unit horizontal vectors $Y$,  (i) $\Ric (Y,U) = 0$ and (ii)  $\Ric (Y,Y) = \Ric (U,U)$.  We now employ the expressions for the Ricci curvature given by Lemmata \ref{lem:vert}, \ref{lem:mix} and \ref{lem:hor}, as well as the formula for the divergence $\dd^*\Om$ given by 
Lemmata \ref{lem:wedge} and \ref{lem:dstar-omega}.  We omit the details.  \hfill q.e.d.  

\medskip

The following example is instructive, even though the soliton behind it is the trivial one $(M^3, g) = \RR^3$.  It is based on an example in \cite{Ba-Ea}.

\begin{example}  {\rm (Helix example)  Choose cylindrical coordinates for $\RR^3$: $(r, \al , z)$, where $r^2 = x_1{}^2 + x_2{}^2,\, \tan \al = x_2/x_1$ and $z = x_3$.  Define $\phi : \RR^3 \ra \RR^2$ by
$$
\phi (r, \al , z) = \left( \ln\left\{ \frac{\sqrt{1+cr^2} - 1}{\sqrt{c}\,r}\right\} + \sqrt{1+cr^2}, - \al + \sqrt{c}\,z\right)\,,
$$
where $c$ is an arbitrary positive constant.  Then $\phi$ is semi-conformal with fibres helices which wind around the concentric cylinders $r=$ constant.  Its gradient is given by
$$
\grad \phi = \left( \frac{\sqrt{1+cr^2}}{r}\frac{\pa}{\pa r}, - \frac{1}{r^2}\frac{\pa}{\pa \al} + \sqrt{c}\frac{\pa}{\pa z}\right)\,,
$$
and its dilation by $\la^2 = (1+cr^2)/r^2$.  Let $(u,v)$ denote coordinates on the codomain $\RR^2$; then $u = u(r)$ with $\ds \frac{\dd u}{\dd r} = \frac{\sqrt{1+cr^2}}{r}$, whereas $v = v(\al , z)$.  The objects associated to our construction of solitons are given as follows:
\begin{eqnarray*}
\theta & = & \frac{-\sqrt{c}r^2}{\sqrt{1+cr^2}}\dd\al - \frac{1}{\sqrt{1+cr^2}}\dd z \\
\Om & = & - \frac{\sqrt{c}r(2+cr^2)}{(1+cr^2)^{3/2}}\dd r\wedge \dd\psi + \frac{cr}{(1+cr^2)^{3/2}}\dd r\wedge \dd z \\
\mu^{\flat} & = & \dd \ln \left( \frac{1}{\sqrt{1+cr^2}}\right) \\
\ov{\Om} & = &  = \frac{2\sqrt{c}r^2}{(1+cr^2)^2}\dd u\wedge \dd v\,.
\end{eqnarray*} 
Note that the functions $\la$ and $\rho$ are basic, so we are in the situation of Theorem \ref{th:main}.  It is routine to check that the equations of the Theorem are satisfied.
}
\end{example}

\begin{example} \label{ex:sol} {\rm In this example we generalize a construction of Ivey \cite{Iv-2} to include examples of non-gradient type.  In particular we will see the geometry Sol arising in this way.  Ivey considers doubly warped product metrics of the form $g = \dd x_1{}^2 + a(x_1)^2\dd x_2{}^2 + b(x_1)^2\dd x_3{}^2$ (but where $\dd x_2{}^2$ and $\dd x_3{}^2$ are to be considered as metrics on a sphere and an Einstein manifold, respectively).  Let $N^2 = \RR^2$ with metric $h = \dd x_1{}^2 + a(x_1)^2 \dd x_2{}^2$ and let $\phi : (\RR^3, g) \ra (\RR^2, h)$ be the projection $\phi (x_1, x_2, x_3) = (x_1, x_2)$.  Then $\phi$ is semi-conformal with dilation $\la \equiv 1$ and with corresponding form dual to its kernel given by $\theta = b(x_1)\dd x_3$.  The dual of the mean curvature of the fibres is given by the form $\mu^{\flat} = \dd \ln b^{-1}$ so that $\rho = b^{-1}$.  Then $\wt{\Om} = \dd \theta + \dd\ln\rho\wedge \theta = 0$ and $\Hh$ is integrable.  In fact Ivey takes the vector field $E$ in (\ref{eq:soliton}) to be the gradient of a function which depends on $x_1$ only.  This would correspond to the vanishing of the function $f$ in the decomposition $E = X + fU$. From (\ref{eq:parts2-bis}), we see that this corresponds to the vanishing of the constant on the right-hand side of equation (ii)(b) of Theorem \ref{th:main}.  We will suppose that the function $\ov{\nu} = \ov{\nu}(x_1)$, but we could conceivably have $f$ non-zero, leading to a more general situation.  On noting that $K^N = - a^{\prime\prime}/a$ and writing $\be = \ln\ov{\nu}$, the equations for a soliton become:
$$
\left\{ \begin{array}{lrl}
{\rm (i)} & - \frac{a^{\prime\prime}}{a} + \frac{1}{2}\left( \be^{\prime\prime}+\frac{a^{\prime}}{a}\be^{\prime} - \left(\frac{b^{\prime}}{b}\right)^2\right) + A & = 0 \\
{\rm (ii)(b)} &   \left(\frac{b^{\prime}}{b}\right)^2 - \frac{b^{\prime\prime}}{b} - \frac{a^{\prime}b^{\prime}}{ab} + \frac{b^{\prime}}{b}\be^{\prime} & = {\rm const.} \\
{\rm (u)} & \be^{\prime\prime} - \left(\frac{b^{\prime}}{b}\right)^2  & = \frac{a^{\prime}}{a} \be^{\prime}\,.
\end{array} \right.
$$

On setting $a = b$, we retrieve the solitons of Example \ref{ex:wp}.
Another particular solution is given by $a = b^{-1}$ with $a = e^{x_1}$, $A = 2$ and the constant on the right-hand side of (ii)(b) equal to $2$.  This gives the metric $g = \dd x_1{}^2 + e^{2x_1}\dd x_2{}^2 + e^{-2x_1}\dd x_3{}^2$, which corresponds to the $3$-dimensional geometry Sol. We can be more explicit in describing the vector field $E$ and seeing if the soliton is of gradient type.  For this we apply equation (\ref{eq:soliton-parts})(ii), which takes the form
$$
\dd f + f\,\dd x_1 + 4e^{-x_1}\dd x_3 = 0\,.
$$
This is solved by
$$
f = - 4x_3e^{-x_1}\,.
$$
Now $E = - \mu + fU + \grad \ln \nu$, which is (locally) a gradient if and only if $\dd E^{\flat} = 0$, i.e. if and only if $\dd (f\theta )= \dd f \wedge \theta + f\Om = 0$.   But it is easily checked that $\dd f \wedge \theta + f\Om = 8x_3e^{-2x_1}\dd x_1\wedge \dd x_3$, which is clearly non-vanishing, so Sol, viewed as a soliton in this way, is not of gradient type.  In fact we can exhibit its soliton flow:
\begin{equation} \label{eq:sol-flow}
E = - 4x_3\frac{\pa}{\pa x_3} - 2\frac{\pa}{\pa x_1}\,.
\end{equation} 
}
\end{example}

\section{The case of minimal fibres} \label{sec:minimal}  

When the potential function $\rho$ is constant, the fibres of $\phi$ are minimal and the equations for a soliton become more accessible; in this case, we are able to give a complete local description of the solutions as specified by Corollary \ref{cor:main}.  Suppose then that $\rho$ is constant.  From (\ref{eq:main})(ii)(b), we see that $\ov{\psi}$ is also a constant,
 $C$ say.  The equations (\ref{eq:main}) now become 
\begin{equation}\label{eq:main-bis}
\left\{ \begin{array}{rrcl}
{\rm (i)} & \Delta^N\ln\ov{\la} + K^N + \frac{1}{2}\Delta^N\ln\ov{\nu} + \frac{(A-C)}{\ov{\la}^2} & = & 0 \\
{\rm (ii)} & \Delta^N\ln\ov{\nu} + \frac{(C+A)}{\ov{\la}^2} & = & 0 \\
{\rm (u)} & \na \dd \ln\ov{\nu} + 2\dd\ln\ov{\la} \odot \dd\ln\ov{\nu}& = &\al \, h\,,
\end{array} \right.
\end{equation}
for some function $\al : N \ra \RR$, with (ii) vacuous whenever $\Hh$ is integrable, i.e. $C=0$.

\begin{proposition} \label{prop:sol}  
Let $\ov{\la}, \ov{\nu}$ satisfy the system {\rm (\ref{eq:main-bis})} with $C\neq 0$.  Then either $C+A=0$ and the corresponding soliton has constant curvature, or $3C-A=0$ and, in terms of a local isothermal coordinate $z= x+\ii y$ on a simply connected domain of $N^2$, we have
\begin{equation} \label{eq:min}
\ov{\la}(z) = \frac{B}{|v(z)|}, \qquad \frac{\pa \ln\ov{\nu}}{\pa z} = \frac{Cv(z)}{B^2} \int\ov{v(z)}\dd \ov{z}\,,
\end{equation}
where $v(z)$ is a non-vanishing holomorphic function and $B$ is a positive constant.
\end{proposition}

\noindent \emph{Remark}:  In fact equations (\ref{eq:main}) and (\ref{eq:main-bis}) depend only on the gradient of $\ga = \ln\ov{\nu}$, so we only expect to be able to explicitly express $\ga_x - \ii \ga_y = 2\frac{\pa\ga}{\pa z}$, as in the proposition above.

\medskip

\noindent \emph{Proof}:  First note that, on combining Proposition \ref{prop:cc} with equation (\ref{eq:main-bis})(i), we see that a solution has constant curvature if and only if $C+A=0$.  In what follows, write $\be = \ln\ov{\la}, \ga = \ln\ov{\nu}$.  Since solutions to the system (\ref{eq:main-bis}) are invariant under conformal changes $h\ra \wt{h} = e^{2u}h$, and since we can always choose a local isothermal coordinate $z = x+\ii y$ with respect to which $h = \delta (z)^2(\dd x^2 + \dd y^2)$, it is no loss of generality to work on a simply connected domain $U \subset \RR^2$ with metric $h = \dd x^2 + \dd y^2$.  Write $\ga_z = \frac{\pa\ga}{\pa z} = \frac{1}{2}\left(\frac{\pa \ga}{\pa x} - \ii \frac{\pa \ga}{\pa y}\right)$, etc.. The system (\ref{eq:main-bis}) now takes the form
\begin{equation} \label{eq:main-complex}
\left\{ \begin{array}{lrl}
{\rm (i)} & 4\be_{z\ov{z}} - \left( \frac{3C-A}{2}\right) e^{-2\be} & = 0 \\
{\rm (ii)} & 4\ga_{z\ov{z}} + (C+A)e^{-2\be} & = 0 \\
{\rm (u)} & \ga_{zz} + 2\ga_z\be_z & = 0 \,.
\end{array}\right.
\end{equation}
On differentiating (u) with respect to $\ov{z}$, (ii) with respect to $z$ and combining, we obtain
$$
\ga_z \be_{z\ov{z}} = 0\,.
$$
Thus either $\ga_z \equiv 0$ on some open set, in which case from (ii) we have $C+A=0$, so that the corresponding soliton has constant curvature, or $\be_{z\ov{z}} \equiv 0$ on some open set, which from (i) implies that $3C-A=0$.  Henceforth we will suppose that we are in the latter situation, so that $\be_{z\ov{z}}\equiv 0$ and $A=3C$.  Note that since $C\geq 0$, this implies that $A\geq 0$ and any corresponding soliton is expanding.  The system (\ref{eq:main-complex}) now takes the form
\begin{equation} \label{eq:main-complex-spe}
\left\{ \begin{array}{lrl}
{\rm (i)} & 4\be_{z\ov{z}} & = 0 \\
{\rm (ii)} & \ga_{z\ov{z}}  & = Ce^{-2\be} \\
{\rm (u)} & \ga_{zz} + 2\ga_z\be_z & = 0 \,.
\end{array}\right.
\end{equation}

If $\ga_z \equiv 0$, then $C=0$, which in turn implies that $A=0$ and we are in the constant curvature case again.  So we will work on a domain where $\ga_{z}\neq 0$.  Then
\begin{equation} \label{eq:gamma-z}
\frac{\ga_{zz}}{\ga_z} = - 2\be_z \quad  \Rightarrow \quad \frac{\pa}{\pa z}\ln \ga_z = - 2\be_z \quad \Rightarrow \quad \ga_z = \mu (z,\ov{z})e^{-2\be}\,,
\end{equation}
where $\frac{\pa \mu}{\pa z} = 0$, i.e. where $\mu = \mu (\ov{z})$ is \emph{antiholomorphic} in $z$.  Then 
$$
\ga_{z \ov{z}} = \left( \mu^{\prime} - 2\mu \be_{\ov{z}}\right)e^{-2\be}\,,
$$
and from (\ref{eq:main-complex-spe}) we require that $C = \mu^{\prime} - 2\mu \be_{\ov{z}}$.  Thus
$$
\be_{\ov{z}} = \frac{\mu^{\prime} - C}{2\mu}\,,
$$
with $\mu = \mu (\ov{z})$ antiholomorphic.  Since the right-hand side of this expression is antiholomorphic, we can integrate along any path in a simply connected domain.  To write this integral more conveniently, we let $p(\ov{z})$ be an antiholomorphic function which is a primitive of $1/\mu$, i.e. $p^{\prime}(\ov{z}) = 1/\mu (\ov{z})$.  Then
$$
\be = \int \frac{\mu^{\prime} - C}{2\mu}\cdot \dd \ov{z} =  - \frac{1}{2}\ln p^{\prime}(\ov{z}) - \frac{C}{2}p(\ov{z}) + q(z)\,,
$$
with $q(z)$ holomorphic in $z$ ($\frac{\pa q}{\pa \ov{z}} = 0$).

Let us change notation once more, by first setting $r(\ov{z}) = e^{p(\ov{z})},\, s(z) = e^{q(z)}$.  Then
$$
\be = - \frac{1}{2} \ln \left\{ \left(\frac{r^C}{C}\right)^{\prime}s^{-2}\right\}\,.
$$
Finally, write $u(\ov{z}) = \left( \frac{r^C}{C}\right)^{\prime}$ and $v(z) = s(z)^{-2}$, to obtain
$$
\be = - \frac{1}{2}\ln \left( u(\ov{z})v(z)\right)\,,
$$
where we take the principal branch of $\ln$ and require that $u(\ov{z})v(z)$ be real and positive for all $z$.  But this is easily seen to be equivalent to the condition
$$
u(\ov{z}) = a \ov{v(z)}\,,
$$
with $a$ a real positive constant.  We now have
$$
\be = - \frac{1}{2}\ln \left( a|v(z)|^2\right)\,,
$$
where $v(z)$ is a holomorphic function, which we require to be non-vanishing to avoid singular points.  In particular, this gives
$$
\ov{\la} = \frac{B}{|v(z)|}\,,
$$
with $B$ a positive constant.

Returning to $\ga_z$; from (\ref{eq:gamma-z}), this is given by
$$
\ga_z = \mu e^{-2\be} = \frac{1}{p^{\prime}}a|v(z)|^2\,.
$$
But $p^{\prime} = \frac{r^{\prime}}{r} = \frac{u}{r^C}$, with $r^C = \int Cu(\ov{z})\dd \ov{z}$.  On noting that $a^{-1/2} = B$, we obtain 
$$
\ga_z = \frac{Cv(z)}{B^2}\int \ov{v(z)}\, \dd \ov{z}\,,
$$
as required.  \hfill q.e.d.    

\begin{corollary} \label{cor:sol}  Any solution to the equations {\rm (\ref{eq:min})} determines a Riemannian $3$-manifold which is locally isometric to the geometry {\rm Nil}.
\end{corollary} 

\noindent \emph{Proof} :  We apply the ansatz of Theorem \ref{th:main}.  Without loss of generality, we may take the constants $B = 1$ and $C = \frac{1}{2}$.  Then
$$
\ov{\Om} = \frac{\ii}{2} v(z) \ov{v(z)} \, \dd z \wedge \dd \ov{z} = \dd \left( \frac{\ii}{2} \left(\int v(z)\, \dd z \right) \ov{v(z)}\, \dd \ov{z}\right)\,.
$$
Let $u(z)$ be a primitive for $v(z)$ : $u^{\prime}(z) = v(z)$.  Then $\ov{\Om} = \dd \left( \Re \left\{ \frac{\ii}{2} u\, \dd \ov{u}\right\} \right)$, and we may set 
$$
\Om = \dd \left( \Re \left\{ \frac{\ii}{2} u \, \dd \ov{u}\right\} + \dd t\right) \,,
$$
to give the metric $g$ in the form
$$
g = \dd u \, \dd \ov{u} + \left( \Re \left\{ \frac{\ii}{2} u \, \dd \ov{u}\right\} + \dd t \right)^2\,.
$$
We claim this is the metric for Nil.  Indeed, if we take $u = y_1 + \ii y_2$ as a local coordinate, then
\begin{equation} \label{eq:gNil}
g = \dd y_1{}^2 + \dd y_2{}^2 + (y_1 \dd y_2 + \dd y_3)^2 
\end{equation}
where we have set $ y_3 = t - \frac{y_1y_2}{2}$.
\hfill q.e.d.  

\bigskip

\noindent \emph{The geometry {\rm Nil} as a soliton}:  On expressing the metric $g$ for Nil as in (\ref{eq:gNil}), 
 we can solve (\ref{eq:soliton-parts})(ii)  for $f$ to obtain the soliton flow $E$ explicitly.  Indeed (\ref{eq:soliton-parts})(ii) becomes
$$ 
\dd f + y_1 \dd y_2 + y_2 \dd y_1 + 2\dd y_3 = 0\,,
$$
which has solution $f =  - y_1y_2 - 2y_3$.  Then the soliton is of gradient type if and only if 
$\dd E^{\flat} = \dd f\wedge \theta + f\Om = 0$.  But
$$
\dd f\wedge\theta + f\Om = - 2(y_1y_2+y_3)\dd y_1\wedge \dd y_2 - y_2\dd y_1\wedge \dd y_3 + y_1\dd y_2\wedge \dd y_3\,,
$$
which is non-vanishing and so Nil, viewed as a soliton in this way, is not of gradient type.  The soliton flow $E$ is given explicitly by
\begin{equation} \label{eq:nil-flow}
E = - y_1\frac{\pa}{\pa x_1} - y_2\frac{\pa}{\pa y_2} - 2y_3\frac{\pa}{\pa y_3} 
\end{equation} 

\medskip

It remains to consider the case when $C = 0$, so that the horizontal distribution is integrable and equations (\ref{eq:main-bis}) become the system:
$$
\left\{ \begin{array}{rrcl}
{\rm (i)} & \Delta^N\ln\ov{\la} + K^N + \frac{1}{2}\Delta^N\ln\ov{\nu} + \frac{A}{\ov{\la}^2} & = & 0 \\
{\rm (u)} & \na \dd \ln\ov{\nu} + 2\dd\ln\ov{\la} \odot \dd\ln\ov{\nu}& = &\al \, h\,,
\end{array} \right.
$$
On replacing $h$ by the conformally related metric $h/\ov{\la}^2$, we may suppose that $\ov{\la} \equiv 1$.  But then the system becomes
$$
\left\{ \begin{array}{rrcl}
{\rm (i)} & K^N + \frac{1}{2}\Delta^N\ln\ov{\nu} + A & = & 0 \\
{\rm (u)} & \na \dd \ln\ov{\nu}  & = &\al \, h\,.
\end{array} \right.
$$
But this is precisely the equation for a $2$-dimensional gradient soliton: $- K^N h = \na\dd \ln\ov{\nu} + Ah$.  This proves Corollary \ref{cor:main}.

\bigskip

\noindent \emph{The case of $\wt{{\rm SL}_2(\RR )}$}:    The geometry $\wt{SL_2(\RR )}$ naturally admits a semi-conformal map with minimal fibres $\phi : \wt{SL_2(\RR )} \ra H^2$ onto the hyperbolic plane.  Specifically, if we write its metric in the form
$$
g = \frac{\dd x_1{}^2 + \dd x_2{}^2}{x_2{}^2} + \left(\frac{\dd x_1}{x_2} + \dd x_3\right)^2\,,
$$
then the projection is given by $\phi (x_1, x_2, x_3) = (x_1, x_2)$, where we take the metric on the codomain to be $h = \frac{\dd x_1{}^2 + \dd x_2{}^2}{x_2{}^2}$.   This fact enables us to establish the following theorem.

\begin{theorem} \label{th:sl2}  The geometry $\wt{{\rm SL}_2(\RR )}$ admits no soliton structure.
\end{theorem}

\noindent \emph{Proof} :  We apply equation (\ref{eq:soliton2}) directly.  We have $\theta = \frac{\dd x_1}{x_2} + \dd x_3$, so that $\Om = \dd \theta =$ $ (\dd x_1 \wedge \dd x_2)/x_2{}^2$ and $\psi = \frac{1}{4} ||I||^2 = \frac{1}{4}||\Om ||^2 = \frac{1}{2}$.  From Lemma \ref{lem:dstar-omega}, $\dd^*\Om = 2\psi \theta = \theta$, so that equation (\ref{eq:soliton2}) becomes
\begin{equation} \label{eq:sl2}
\left( A - \frac{3}{2}\right) g + \frac{1}{2} \Ll_Xg + \dd f \odot \theta + 2 \theta^2 = 0\,.
\end{equation}
Thus $\wt{{\rm SL}_2(\RR )}$ admits a soliton structure if and only if this has a solution for $f$, $X$ and $A$.

Let $Y_1 = x_2 \pa_1 - \pa_3$, $Y_2 = x_2\pa_2$ be an orthonormal frame for $\Hh$, and let $X = \al Y_1 + \be Y_2$.  Then a routine calculation gives
$$
\frac{1}{2}\Ll_X g = \frac{1}{x_2} (\al \dd x_2 - \be \dd x_1)\left( \frac{2\dd x_1}{x_2} + \dd x_3\right) + \frac{\dd \al \dd x_1}{x_2} + \frac{\dd \be \dd x_2}{x_2}\,.
$$
On substituting into (\ref{eq:sl2}) and equating the various coefficients of $\dd x^i \dd x^j$ to zero, we are led to the following system of equations:
\begin{equation} \label{eq:sl2-system}
\left\{ \begin{array}{ll}
{\rm (i)} &  2A - 1 - 2\be + x_2 \pa_1\al + x_2 \pa_1 f = 0  \\
{\rm (ii)} & A - \frac{3}{2} + x_2 \pa_2 \be = 0   \\
{\rm (iii)} & A + \frac{1}{2} + \pa_3 f = 0  \\
{\rm (iv)} & \frac{2\al}{x_2} + \pa_2 \al + \pa_1 \be + \pa_2 f = 0  \\
{\rm (v)} & 2(A + \frac{1}{2}) - \be + \pa_3 \al + \pa_3 f + x_2 \pa_1 f = 0  \\
 {\rm (vi)} & \al + \pa_3\be + x_2 \pa_2 f = 0 
\end{array}
\right.
\end{equation}
We claim this has no solution in $\al , \be , f$, whatever the value of the constant $A$.

It is convenient to write $a = A + \frac{1}{2}$.  Then (iii) and (ii) imply 
\begin{eqnarray*}
f & = & - a x_3 + p(x_1, x_2) \\
\be & = & - (a-2)\ln x_2 + q(x_1, x_3)
\end{eqnarray*}
for functions $p = p(x_1, x_2), q = q(x_1, x_3)$.  From (vi), $\al = - \pa_3q - x_2 \pa_2p$, and we obtain the following system in $p$ and $q$:
\begin{equation} \label{eq:sl2-system2}
\left\{ \begin{array}{ll}
{\rm (vii)} &  2(a-1) + 2(a-2) \ln x_2 - 2q - x_2\pa_{13}q - x_2{}^2\pa_{12}p + x_2 \pa_1p = 0  \\
{\rm (viii)} & - \frac{2}{x_2}\pa_3q - 2\pa_2p - x_2 \pa_{22}p + \pa_1q = 0   \\
{\rm (ix)} & a + (a-2)\ln x_2 - q - \pa_{33}q + x_2 \pa_1p = 0
\end{array}
\right.
\end{equation}
On differentiating (ix) with respect to $x_2$, substituting into (vii) and using (ix) once more, we obtain the equation $a - 4 - x_2 \pa_{13}q + 2\pa_{33}q = 0$,
which, on integrating with respect to $x_3$, implies that
$$
(a - 4)x_3 - x_2 \pa_1q + 2 \pa_3 q = \wt{r}(x_1, x_2)\,,
$$
for some function $\wt{r} = \wt{r}(x_1, x_2)$.  But since $q$ is independent of $x_2$, we must have $- \pa_1q = \pa_2\wt{r}$, which implies that $\wt{r} = - x_2\pa_1q + r(x_1)$, for a function $r = r(x_1)$.  Therefore $(a-4)x_3 + 2\pa_3q = r(x_1)$, which, on integrating with respect to $x_3$, gives
$$
q = \frac{1}{2} \left\{ x_3 r(x_1) + s(x_1) - (a-4)\frac{x_3{}^2}{2}\right\}\,,
$$
for some function $s = s(x_1)$.  But the fact that $\wt{r}$ is independent of $x_3$ shows that $r^{\prime}$ vanishes and so $r$ is constant.  We can now substitute back into the system (\ref{eq:sl2-system2}) to obtain:
\begin{equation} \label{eq:sl2-system3}
\left\{ \begin{array}{ll}
{\rm (x)} &  2(a-1) + 2(a-2) \ln x_2 - x_3r - s + \frac{1}{2}(a-4)x_3{}^2 - x_2{}^2 \pa_{12}p + x_2 \pa_1p = 0  \\
{\rm (xi)} & - \frac{1}{x_2}\left( r - (a-4)x_3\right) - 2\pa_2p - x_2 \pa_{22}p + \frac{1}{2} s^{\prime} = 0   \\
{\rm (xii)} & a + (a-2)\ln x_2 - \frac{1}{2}\left( rx_3 + s - \frac{a-4}{2} x_3\right) + \frac{1}{2}(a-4) + x_2 \pa_1 p = 0
\end{array}
\right.
\end{equation}
Then (xi) implies that $a = 4$, which, on substituting into (xii) gives $r = 0$.  From (xi), we have $2\pa_2p + x_2\pa_{22}p = \frac{1}{2}s^{\prime}$, which implies that $x_2\pa_2p + p = \frac{x_2}{2}s^{\prime} + t(x_1)$, for some function $t = t(x_1)$.  But then $x_2 \pa_{12}p + \pa_1p = \frac{x_2}{2}s^{\prime\prime} + t^{\prime}$.  Equations (x) and (xii) now quickly lead to a contradiction.
\hfill q.e.d. 

\bigskip

\noindent \emph{The other geometries}:  We can apply the same techniques as above to the other geometries (and in principle to any $3$-manifold admitting a semi-conformal map to a surface), to obtain a system of equations similar to (\ref{eq:sl2-system}). The calculations are long and tedious so we omit the details, but a complete set of solutions can be obtained.  We summarize the conclusions as follows:  \emph{up to addition of a Killing vector field, the soliton flows $E$ on {\rm Nil} given by {\rm (\ref{eq:nil-flow})} and on {\rm Sol} given by {\rm (\ref{eq:sol-flow})}, are unique;  the soliton flows on $S^2 \times \RR$ and $H^2 \times \RR$ given in {\rm Example \ref{ex:wp}} are unique; the only soliton flows on $S^3$ and $H^3$ are given by Killing vector fields; the only non-Killing solitons on $\RR^3$ are the well-know Gaussian solitons}.

\end{document}